\documentclass[11pt]{article}
\usepackage{amscd}

\let\ssection=\section
\renewcommand{\section}{\setcounter{equation}{0}\ssection}

\newcommand{\Diff}{\mathrm{Diff}}

\newcommand{\cF}{{\mathcal{F}}}
\newcommand{\cD}{{\mathcal{D}}}

\newcommand{\Pol}{\mathrm{Pol}}

\newcommand{\SL}{\mathrm{SL}}
\newcommand{\Sl}{\mathrm{sl}}

\newcommand{\Vect}{\mathrm{Vect}}

\def\d{\delta}
\def\G{\Gamma}
\def\om{\omega}

\def\a{\alpha}
\def\b{\beta}

\def\l{\lambda}
\def\m{\mu}
\def\n{\nabla}
\def\p{\partial}

\newtheorem{thm}{Theorem}[section]

\newtheorem{cor}[thm]{Corollary}
\newtheorem{pro}[thm]{Proposition}

\newtheorem{rmk}[thm]{Remark}

\newcommand{\bbR}{\mathbf{R}}

\newcommand{\bbC}{\mathbf{C}}

\begin{document}

\title{Projectively Equivariant Quantization Map }

\author{
{\large Sofiane BOUARROUDJ}\\
{\footnotesize {\it  CNRS, Centre de Physique Th\'eorique, Luminy, Case 907,
F13288 Marseille, Cedex 9, France.}}\\
{\footnotesize { \it e-mail: sofbou@cpt.univ-mrs.fr}}
}
\date{}
\maketitle

\begin{abstract}
 Let $M$ be a manifold endowed with a symmetric affine connection
$\G.$ The aim of this Letter is to describe a quantization map between the
space of second-order polynomials on the cotangent bundle $T^{*} M$ and the
space of second-order linear differential operators, both viewed as modules
over   the group of diffeomorphisms and the Lie algebra of vector fields on
$M.$  This map is  an isomorphism, for almost all values of certain
constants, and it depends only on the projective class of the
affine connection~$\G.$
\end{abstract}

\section{Introduction}
The quantization procedure proposed in this Letter deals with
the space of linear differential operators and the corresponding
space of symbols viewed as modules over the group of diffeomorphisms
$\Diff(M)$ and the Lie algebra of vectors fields $\Vect(M).$
This method of quantization have been introduced in
the recent papers (\cite{lo}, \cite{dlo}, \cite{do2}).

Let $\cD_{\l,\mu}(M)$ be the space of linear differential operators
from the space of $\l-$densities to the space of $\mu-$densities. The
corresponding space of symbols, $\Pol_{\d}(T^*M),$ is the space of
fiberwise polynomials on $T^*M$  with
values in the space of $\d-$densities over $M$, where~$\d=\mu-\l.$

We call quantization map, a linear map
\begin{equation}
\label{sah}
Q_{\l,\mu}:\Pol_{\d}(T^*M)\rightarrow \cD_{\l,\mu}(M),
\end{equation}
that is bijective and preserves the principal symbol (see \cite{lo},
\cite{do2}, \cite{dlo}).

There is no quantization map (\ref{sah})
equivariant with respect to the action of the
group $\Diff(M)$.
It is natural to consider a subgroup $G\subset \Diff(M)$ (of finite
dimension) and to restrict the action of $\Diff(M)$ on the subgroup $G.$
There are two interesting cases:

If $M=\bbR^n$ is endowed with a flat projective structure,
the quantization map is given in \cite{lo}. This map is equivariant with
 respect to the action of the group of projective transformations
 $\SL_{n+1}\subset \Diff(\bbR^n)$.
If $M=\bbR^n$ is endowed with a flat conformal structure, the
quantization map is given in \cite{do2}, it is equivariant with
 respect to the action of the group of conformal diffeomorphisms
$SO(p+1,q+1)\subset \Diff(\bbR^n),$ where $(p+q=n).$ (See also \cite{ga},
\cite{bo1}, \cite{cmz}, for the one dimensional case.)

A natural and well-known way to define a quantization map is to fix an affine
connection on $M$ (see, e.g., \cite{bgv}). However, there is no canonical
quantization map associated to the given connection.

The purpose of this Letter is to study the quantization map (\ref{sah})
between the space of symbols of degree less than two and the space of
second-order linear differential operators satisfying the following
properties:

1. It is projectively invariant, i.e. it  depends only on the projective
class of the affine connection $\G.$

2. If $M=\bbR^n$ with a flat projective structure, this isomorphism is
equivariant with respect to the action of the group of projective
transformations $\SL_{n+1}$ (resp. infinitesimal projective transformations
$\Sl_{n+1}$).

The method used in this paper follows that of the recent preprint \cite{do}.
\section{Space of Linear Differential Operators}
Let $M$ be a manifold of dimension $n$ endowed with an affine
connection $\Gamma.$
We are interested in defining a two parameter family of $\Diff(M)-$module
(resp. $\Vect(M)-$module) on the space of linear differential operators. This
space was recently studied in recent papers
(\cite{bo1}, \cite{bo2}, \cite{ga}, \cite{lo}, \cite{dlo}, \cite{do},
\cite{do2}, {\cite{lmt}).
\subsection{Space of Tensor Densities}
For simplicity, we assume that $M$ is  oriented throughout this Letter.

The space of tensor densities on $M,$  $\cF_{\l}(M),$ or $\cF_{\l}$ for
simplify, is the space of sections of the line bundle
$(\Lambda^n T^*M)^{\otimes \lambda},$ where $\l \in \bbC.$ As a vector space,
tensor densities are isomorphic to the space of complexified
functions, but the structure
of $\Diff(M)$-module is different. Let us explicate this action:

Let $f\in \Diff(M)$ and  $\phi \in
{\cal F}_{\l}.$ In a local coordinates $(x^i)$,  the action is given by
\begin{equation}
\label{den}
f^*\phi=\phi\circ f^{-1}\cdot ( {J_{f^{-1}}})^{\lambda},
\end{equation}
where $J_f=\left |\frac{Df}{Dx}\right |$ is the Jacobian of $f$.\\
In the case $\l =0,$ $1$, the action (\ref{den}) is precisely
the standard action of $\Diff(M)$ on the space of functions and differential
forms of degree $n$ respectively.

Differentiating the action of the flow of a vector field,  one gets
the corresponding  represen\-tation of $\Vect(M).$
\begin{equation}
\label{denv}
L^{\lambda}_X(\phi)=X^i\partial_i(\phi)+\lambda\partial_i(X^i)\phi,
\end{equation}
where $X=X^i\p_i\in \Vect(M).$

Formul\ae \,(\ref{den}), (\ref{denv}) do not depend on the choice of
coordinates.\\

Let us now recall the definition of covariant derivative on tensor densities
(cf. \cite{do2}).\\
Let $\n$ be the covariant derivative associated to the affine connection
$\Gamma.$
  If $\phi \in \cF_{\l}$, then $\n \phi \in \Omega^{1}(M)\otimes \cF_{\l},$
given, in a local coordinates, by the formula
\begin{equation}
\label{covden}
\n_i\phi=\p_i \phi-\l \G_i\phi,
\end{equation}
with $\G_i=\G_{il}^l$ (summation is understood on repeated indices).
\subsection{Space of Linear Differential Operators}
Consider $\cD_{\l, \mu}(M)$, the space of linear differential operators
acting on tensor densities
\begin{equation}
A:\cF_\l\to\cF_\m.\nonumber
\label{Conv}
\end{equation}
The action of $\Diff(M)$ on $\cD_{\l, \mu}(M)$ depends on two parameters
$\l$ and $\m$.
This action is given by the equation:
\begin{equation}
f_{\l,\m}(A)=f^*\circ A\circ {f^*}^{-1},
\label{Opaction}
\end{equation}
where $f^*$ is the action (\ref{den}) of $\Diff(M)$ on $\cF_\l$.

Differentiating the action of the flow of a vector field,  one gets
the corresponding  represen\-tation of $\Vect(M).$
\begin{eqnarray}
\label{noir et blanc}
L^{\l,\mu}_X(A)=L^{\mu}_X\circ A-A\circ L^{\lambda}_X,
\end{eqnarray}
where $X\in \Vect(M).$

These formul\ae\, do not depend on the choice of a system of coordinates.

Denote ${\cal D}^k_{\l,\mu}(M)$ the space of $k-$order
linear differential operators acting on tensors densities. In a local
coordinates $(x^i),$ one can write  $A\in {\cal D}_{\l, \mu}^k(M)$
\begin{equation}
\label{un truc}
A=a_k^{i_1,\ldots,i_k} \frac{\partial}{\partial x^{i_1}}\cdots \frac{\partial}
{\partial x^{i_k}}+\cdots +a_1^i \frac{\partial}{\partial x^i}+a_0,\nonumber
\end{equation}
with the coefficients $a_l^{i_1,\ldots,i_l}=a_l^{i_1,\ldots,i_l}
(x^1,\ldots,x^n)\in
C^{\infty}(M),$ with $l=0,\ldots,k.$ We have then a filtration
$$
{\cal D}^0_{\l,\mu}(M)\subset {\cal D}^1_{\l,\mu}(M)\subset \cdots \subset
{\cal D}^k_{\l,\mu}(M)\subset \cdots $$

The space of $k-$order linear differential operators $\cD_{\l,\m}^k(M)$ is a
$\Diff(M)$-submodule (resp. $\Vect(M)-$submodule) of $\cD_{\l,\m}(M).$
\begin{rmk}{\rm The space of linear differential operators viewed as a
module over the group
of diffeomorphisms is a classical object (see e.g., \cite{wil}). For example,
in the case $M=S^1,$ the space of Sturm-Liouville operators
$\frac{d^2}{dx^2}+u(x),$ where $u(x)\in \cF_{2}$, is viewed as a submodule of
$\cD^2_{-\frac{1}{2},\frac{3}{2}}(S^1).$
Also, the modules $\cD^k_{\frac{1-k}{2},\frac{1+k}{2}}(S^1)$ was considered in
\cite{wil}.

}
\end{rmk}
\section{Space of Symbols}
The space of symbols, $\Pol (T^*M),$ is the space of
functions on the cotangent bundle $T^*M$ polynomial on the
fibers. In a local coordinate system $(x_i,\xi_i),$ one can write
$$
T=\sum_{l=0}^kT^{i_1,\ldots,i_l}\xi_{i_1}\cdots \xi_{i_l},
$$
with $T^{i_1,\ldots,i_l}(x^1,\ldots,x^n) \in C^{\infty}(M).$

One defines a one parameter family of $\Diff (M)-$module
(resp. $\Vect(M)-$module) on the space of symbols by
$$
\Pol_{\d}(T^*M):=\Pol(T^*M)\otimes \cF_{\d}.
$$
Let us explicit this action.

Take  $f\in \Diff(M),$ and $X\in \Vect(M)$.
Then, in a local coordinates  $(x_i,\xi_i),$ one has:
\begin{eqnarray}
f_{\d}(T)&=& f^*T\cdot (J_{f^{-1}})^{\d},\\
L^{\d}_X(T)&=& L_X(T)+\d \mathrm{D}(X)\,T, \quad  T\in \Pol_{\d}(T^*M)
\end{eqnarray}
where $$L_X=X^i\p_i-\xi_j \p_i(X^j)\p_{\xi_i},
\quad \mathrm{D}(X)=\p_iX^i.$$

The space of symbols admits a graduation
$$\Pol_{\d}(T^*M)=\bigoplus_{k=0}^{\infty} \Pol_{\d,k}(T^*M),$$
where $ \Pol_{\d,k}^k(T^*M)$ are the homogeneous polynomials of degree $k$ on
$T^*(M).$ This graduation is $\Diff(M)-$invariant.

Throughout this Letter, we will identify  the space of symbols
with the space of symmetric contravariant tensor fields on $M.$
\section{Flat Projective Structure and Projectively Equivalent Connection}
Let $M$ be a manifold of dimension $n.$ Recall two notions on projective
geometry, the notion of flat projective structure and the notion of
projectively equivalent connections  (see \cite{kn}).
\subsection{Flat projective structure}
A manifold $M$ admits {\it a flat projective structure} if there exists
an atlas $\{\phi_i\}$ such that the local transformations
$\phi_i\circ \phi^{-1}_j$ are projective transformations.

The most interesting case is when $M=\bbR^n.$ In this case, the
group $\SL_{n+1}$ acts locally on $\bbR^n$ by projective transformations.
Choosing a local coordinate system, the Lie algebra $\Sl_{n+1}$ can be
identified with the subalgebra of $\Vect(\bbR^n)$ generated by the vector
fields:
$$
\p_i, \quad x^i\p_j,\quad x^ix^j\p_j.
$$

The projective Lie algebra $\Sl_{n+1}$ is a maximal subalgebra of the
Lie algebra of polynomial vector fields on $\bbR^n$ (cf. \cite{lo}).
\subsection{Projectively Equivalent Connections}
\label{connexion}
The notion of projectively equivalent connection is an old notion related to
projective geometry of the
``paths'' studied by H. Weyl in \cite{Weyl} and T.Y. Thomas in \cite{Thomas}.
Weyl gives the following definition:\\
Two affine connection without torsion, with Christoffel symbols
$\tilde{\Gamma}_{jk}^i$ and
$\Gamma_{jk}^i$ given on the same system of coordinate $x^1,\ldots,x^n$, are
projectively equivalent, if there exists a differential 1-form with
components $\om_i,$ such that
\begin{equation}
\label{con}
\tilde{\Gamma}_{jk}^i=\Gamma_{jk}^i+\delta^i_j\om_k+\delta^i_k\om_j.
\end{equation}

Geometrically, two affine connections without torsion projectively equivalent
give the same unparameterized geodesics (cf. \cite{kn}, \cite{Weyl}).

An affine connection $\G$ is said to be {\it projectively flat}, if they can
be written:
\begin{equation}
\label{fla}
\G^i_{jk}=\frac{1}{n+1}\left ( \d^i_j\G_k+ \d^i_k\G_j \right ).
\end{equation}

A manifold $M$ endowed with an affine connection $\G$, admits a flat
projective structure if and only if the connection $\G$ is projectively
flat (cf. \cite{kn}).
\section{Main Theorems}
In this section, we will give a quantization map between the space of
 symbols of degree less than two and the space of second-order linear
differential operators. First, decompose the space of symbols into a direct
sum
$$\Pol_{\d,2}(T^*M)\oplus \Pol_{\d,1}(T^*M)\oplus \Pol_{\d,0}(T^*M) ,$$
where $\Pol_{\d,k}(T^*M)$, is the space of homogeneous
symbols of degree k, ( $k=0,1,2).$  We will construct a quantization map
on each of these spaces.

In the case of symbols of degree less than 1, there exists a quantization map
that commutes with the action of $\Diff(M)$ and $\Vect(M)$
(see \cite{lo}, \cite{do}).
\begin{thm}
\label{th1}
For any $\d \not=1,$ the map
$Q_{\l,\mu}^1: \Pol_{\d, 1}(M)\oplus\Pol_{\d, 0}(M) \rightarrow
\cD_{\l,\mu}^{1}(M)$
given by
\begin{equation}
\label{mainun}
Q_{\l,\mu}^1(T)=T^i\,\n_i+\a\,\n_i(T^i)+T^0
\end{equation}
where $T=T^i\,\xi_i+T^0,$ and
\begin{equation}
\label{acon}
\a= \frac{\l}{1-\d}
\end{equation}
is a projectively  invariant isomorphism (i.e. it depends only on the
projective class of the affine connection $\Gamma.$)
\end{thm}
{\bf Proof of Theorem \ref{th1}.}
Let $\tilde \Gamma$ be a symmetric affine connection projectively equivalent
to $\Gamma.$ Denote by $\tilde Q^1_{\l,\mu}$ the quantization map written
with the connection $\tilde \G$. We must show
$$\tilde Q^1_{\l,\mu}=Q^1_{\l,\mu}.$$

We need some formul\ae \, (see \cite{do2}):\\
Recall the covariant derivative on the space of  1-order contravariant tensor
fields: Let $T^i$ be a tensor, then one has:
\begin{eqnarray}
\label {covsym}
\n_j(T^i)&=&\p_j\, T^i+\Gamma_{jl}^iT^l-\d\, \G_j T^i.
\end{eqnarray}
Using  formul\ae\, (\ref{con}), (\ref{covsym}) one obtains
$$\tilde \n_i(\phi)=\n_i(\phi)-\l(n+1)\,\om_i\,\phi, \quad
\tilde \n_i(T^i)=\n_i(T^i)+(1-\d)(1+n)\,\om_i\, T^i.$$

Then after a straightforward calculation:
$$\tilde Q^1_{\l,\mu}(T)=Q^1_{\l,\mu}(T)+(1+n)(\alpha (1-\d)-\l)\,\om_i\, T^i.$$
Hence $\tilde Q^1_{\l,\mu}=Q^1_{\l,\mu}$ if and only if $\alpha$ is given as
in (\ref{acon}).\\
\begin{rmk}
\label{cor1} {\rm
1. If $M$ is endowed with a flat projective structure, the isomorphism
(\ref{mainun}) is the unique provided it preserves the principal symbol
(cf. \cite{lo}).\\
2. In the case $\d=1,$ the modules are still isomorphic if
$(\l,\mu)=(0,1).$ This isomorphism is given by (\ref{mainun}) with an
arbitrary $\alpha.$\\
}
\end{rmk}
Let us give the quantization map on the space of homogeneous symbols of
degree $2.$
\begin{thm}
\label{cafe} If $n\geq 2,$ for any
$\delta\not= \frac{n+3}{n+1}, \frac{n+2}{n+1}, $ there exists a projectively
invariant isomorphism
$Q_{\l, \mu}^2: \Pol_{\delta,2}(T^*M)\rightarrow \cD_{\l, \mu}^2(M)$ given
by
\begin{equation}
\label{main}
Q_{\l, \mu}^{2}(T)=T^{ij}\n_i\n_j+\b_1\n_jT^{ij}\n_i+\b_2\n_i\n_j(T^{ij})+
\b_3R_{ij}T^{ij},
\end{equation}
where $T(\xi)=T^{ij}\xi_i\xi_j,$ the coefficients
$\beta_1,\b_2, \beta_3$ are as follows
\begin{eqnarray}
\label{cons}
\beta_1&=&\frac{2+2\l (n+1)}{2+(1+n)(1-\d)}\nonumber\\
\beta_2&=&\frac{\l(n+1)(1+\l (n+1))}{((1-\d)(1+n)+1)((1-\d)(1+n)+2)}\\[2mm]
\beta_3&=&\frac{\l\,(\mu -1)(n+1)^2}{(1-n)((1-\d)(1+n)+1)}\nonumber
\end{eqnarray}
and $R_{ij}$ denote the components of Ricci tensor of the connection $\Gamma.$

\end{thm}
\begin{cor}
\label{maindeuxp}
If $M$ is endowed with a flat projective structure then:

1. The isomorphism (\ref{main}) has the  form
\begin{eqnarray}
\label{footing}
Q_{\l, \mu}(T)&=&T^{ij}\p_i\p_j+\b_1 \p_j T^{ij}\p_i+\b_2 \p_i\p_jT^{ij},
\end{eqnarray}
where the constants $\b_1,\b_2$ are as in (\ref{cons}).

2. It is the unique map equivariant with respect to the action of $\SL_{n+1}$
(resp. $\Sl_{n+1}$)  that preserves the principal symbols (cf. \cite{lo}).
\end{cor}
{\bf Proof of the Theorem \ref{cafe}} Let $\tilde\Gamma$ be a connection
projectively equivalent to $\Gamma.$ Denote by
$\tilde Q^2_{\l,\mu}$ the  quantization map written with $\tilde \Gamma.$

We need some formul\ae\, (see \cite{do2}):

The covariant derivative on the space of  2-order contravariant
tensor fields reads:
\begin{eqnarray}
\label {covsy}
\n_k(T^{ij})&=&\p_k T^{ij}+\G_{lk}^iT^{lj}+\G_{lk}^jT^{il}-\d \G_k T^{ij}.
\end{eqnarray}
The second-order term in $Q^2_{\l,\mu}$ reads:
\begin{eqnarray}
\label{hmida1}
T^{ij}\p_i\p_j&=&T^{ij}\n_i\n_j+(T^{jk}\Gamma^{i}_{jk}+2\l T^{ij}\Gamma_j)\,
\n_i\\[2mm]
&  &+T^{ij}(\l^2\Gamma_i\G_j+\l\p_i\G_j),\nonumber
\end{eqnarray}
the first-order term in $Q^2_{\l,\mu}$ reads:
\begin{eqnarray}
\label{hmida2}
\p_jT^{ij}\p_i&=&\n_jT^{ij}\n_i-(T^{jk}\G^i_{jk}+(1-\d)T^{ij}\G_j)\,\n_i\\[2mm]
& &+\l (\n_iT^{ij})\G_j-\l T^{ij}\,(\G_{ij}^k \G_k+(1-\d)\G_i\G_j),\nonumber
\end{eqnarray}
and the zero-order part of $Q^2_{\l,\mu}$ reads:
\begin{eqnarray}
\label{hmida3}
\p_j\p_jT^{ij}&=&\n_i\n_jT^{ij}-2(1-\d)(\n_iT^{ij})\G_j-(\n_iT^{jk})\G_{jk}^i\\
& &-T^{ij}(\p_k\G_{ij}^k+(1-\d)\p_i\G_j-2 \G^l_{ik}\G_{jl}^k-(1-2\d)
\G^k_{ij}\G_k-(1-\d)^2\G_i\G_j).\nonumber
\end{eqnarray}

Now after calculation one has:
\begin{eqnarray}
\label{coeff}
\lefteqn{ \tilde Q_{\l,\mu}(T)= Q_{\l,\mu}(T)
+[2\b_2+(1+n)(-\l\b_1+2\eta_\d\b_2)]\n_iT^{ij}\,\om_j} \nonumber\\[2mm]
& & +[2\b_1-2+(1+n)(-2\l +\eta_\d \b_1)]\,\,
T^{ij}\om_j \n_i \nonumber\\[2mm]
&  &+[ (1+n)(-\l+\eta_\d\b_2)+2\b_2+(1-n)\b_3]\,\,
T^{ij}\p_i \,\om_j\\[2mm]
& & +[(1+n)(\l-\eta_\d\b_2)-2\b_2+\b_3(n-1)]\,\,
T^{jk}\,\Gamma^i_{jk}\,
\om_i\nonumber\\[2mm]
&& +[(1+n)^2 (\l^2+\eta_\d (\d\b_2-\l \b_1))+
\!2(1+n)(\l(1\!-\!\b_1)+
\d\b_2)+(n-1)\b_3]\,T^{ij}\om_i\,\om_j\nonumber
\end{eqnarray}
where $\eta_\d =1-\d.$

Hence, $\tilde Q_{\l,\mu}^2=Q_{\l,\mu}^2$ if and only if the constants
 $\b_1,\b_2,\b_3$ are given as in (\ref{cons}).\\

\noindent {\bf Proof of the Corollary \ref{maindeuxp}.} In this case
(see section \ref{connexion}), the connection $\G$ can be written,
in the coordinates of the flat projective structure, in the form
$$\Gamma^{k}_{ij}=\frac{1}{n+1}(\d^k_i\Gamma_j+\d^k_j\G_i).$$
Substituting  this formula in  Equations
(\ref{hmida1}), (\ref{hmida2}), (\ref{hmida3}), and, finally to the map
(\ref{maindeuxp}) one gets the expression (\ref{footing}).

The proof of the part 2) is given in \cite{lo}.

Let us study the particular values of $\d$ called ``resonant'':
\begin{pro}
\label{po}
In the resonant case $\d= \frac{n+2}{n+1}, \frac{n+3}{n+1}$, the
 modules are still isomorphic with the particular values of
$\l, \mu, \b_1,\b_2,\b_3,$ given in {\rm table I} .
\end{pro}
{\bf Proof of the proposition \ref{po}.} Replace the particular values
of $\d$ in the formula (\ref{coeff}). hence,  $Q_{\l, \mu}^2=\tilde
Q_{\l, \mu}^2,$ if and only if the constants $\l, \mu, \b_1, \b_2, \b_3$ is
given as in the table I.
\begin{rmk}{\rm
In contrast with the nonresonant case, if $M$ is flat and
$\d=\frac{n+3}{n+1},$ the isomorphism is not unique. There is a family of
isomorphisms with arbitrary constant $\b_2.$ }
\end{rmk}
\begin{center}
\label{tableau}
\begin{tabular}{|c||c||c||c||c||c|}\hline
$\d$ & $\l$ & $\mu$ \rule[-5mm]{0mm}{12mm}
 & $\b_1$  & $\b_2$  & $\b_3$ \\ \hline
$\displaystyle \frac{n+3}{n+1}$ &
$\displaystyle\frac{-1}{n+1}$ \rule[-5mm]{0mm}{12mm}&
$\displaystyle\frac{n+2}{n+1}$ &
$2\b_2$  & .
  &
$\displaystyle\frac{1}{1-n}$\\[2mm] \hline
$\displaystyle\frac{n+2}{n+1}$ & $0$  &$\displaystyle\frac{n+2}{n+1}$
\rule[-5mm]{0mm}{12mm}   & $2$  &  $0$ & $0$ \\[2mm] \hline

$\displaystyle\frac{n+2}{n+1}$ &   $\displaystyle
\frac{-1}{n+1}$&
\rule[-5mm]{0mm}{12mm}   $1$  &   $0$ &   $0$ &
$\displaystyle\frac{1}{1-n}$
\\[2mm] \hline

\end{tabular}\\

\vspace*{0.5cm}
Table I.
\end{center}

It would be interesting to obtain an analogue of the formula (\ref{main})
in the case of higher-order differential operators.

\bigskip

{\it Acknowledgements}. I am indebted to C. Duval and V. Ovsienko
for the statement of the problem and for numerous fruitful discussions.
\vskip 1cm


\end{document}